\documentclass[reqno,12pt,a4paper]{amsart}

\usepackage{mathrsfs}
\usepackage{amssymb}
\usepackage{amsfonts}
\usepackage{latexsym}
\usepackage{amsthm}
\usepackage{graphicx}
\def\lb{\label}

\newcommand{\er}[1]{\textrm{(\ref{#1})}}

\begin{document}


\renewcommand{\theequation}{\arabic{section}.\arabic{equation}}
\theoremstyle{plain}
\newtheorem{theorem}{\bf Theorem}[section]
\newtheorem{lemma}[theorem]{\bf Lemma}
\newtheorem{corollary}[theorem]{\bf Corollary}
\newtheorem{proposition}[theorem]{\bf Proposition}
\newtheorem{definition}[theorem]{\bf Definition}
\newtheorem{remark}[theorem]{\bf Remark}

\def\a{\alpha}  \def\cA{{\mathcal A}}     \def\bA{{\bf A}}  \def\mA{{\mathscr A}}
\def\b{\beta}   \def\cB{{\mathcal B}}     \def\bB{{\bf B}}  \def\mB{{\mathscr B}}
\def\g{\gamma}  \def\cC{{\mathcal C}}     \def\bC{{\bf C}}  \def\mC{{\mathscr C}}
\def\G{\Gamma}  \def\cD{{\mathcal D}}     \def\bD{{\bf D}}  \def\mD{{\mathscr D}}
\def\d{\delta}  \def\cE{{\mathcal E}}     \def\bE{{\bf E}}  \def\mE{{\mathscr E}}
\def\D{\Delta}  \def\cF{{\mathcal F}}     \def\bF{{\bf F}}  \def\mF{{\mathscr F}}
\def\c{\chi}    \def\cG{{\mathcal G}}     \def\bG{{\bf G}}  \def\mG{{\mathscr G}}
\def\z{\zeta}   \def\cH{{\mathcal H}}     \def\bH{{\bf H}}  \def\mH{{\mathscr H}}
\def\e{\eta}    \def\cI{{\mathcal I}}     \def\bI{{\bf I}}  \def\mI{{\mathscr I}}
\def\p{\psi}    \def\cJ{{\mathcal J}}     \def\bJ{{\bf J}}  \def\mJ{{\mathscr J}}
\def\vT{\Theta} \def\cK{{\mathcal K}}     \def\bK{{\bf K}}  \def\mK{{\mathscr K}}
\def\k{\kappa}  \def\cL{{\mathcal L}}     \def\bL{{\bf L}}  \def\mL{{\mathscr L}}
\def\l{\lambda} \def\cM{{\mathcal M}}     \def\bM{{\bf M}}  \def\mM{{\mathscr M}}
\def\L{\Lambda} \def\cN{{\mathcal N}}     \def\bN{{\bf N}}  \def\mN{{\mathscr N}}
\def\m{\mu}     \def\cO{{\mathcal O}}     \def\bO{{\bf O}}  \def\mO{{\mathscr O}}
\def\n{\nu}     \def\cP{{\mathcal P}}     \def\bP{{\bf P}}  \def\mP{{\mathscr P}}
\def\r{\rho}    \def\cQ{{\mathcal Q}}     \def\bQ{{\bf Q}}  \def\mQ{{\mathscr Q}}
\def\s{\sigma}  \def\cR{{\mathcal R}}     \def\bR{{\bf R}}  \def\mR{{\mathscr R}}
\def\S{\Sigma}  \def\cS{{\mathcal S}}     \def\bS{{\bf S}}  \def\mS{{\mathscr S}}
\def\t{\tau}    \def\cT{{\mathcal T}}     \def\bT{{\bf T}}  \def\mT{{\mathscr T}}
\def\f{\phi}    \def\cU{{\mathcal U}}     \def\bU{{\bf U}}  \def\mU{{\mathscr U}}
\def\F{\Phi}    \def\cV{{\mathcal V}}     \def\bV{{\bf V}}  \def\mV{{\mathscr V}}
\def\P{\Psi}    \def\cW{{\mathcal W}}     \def\bW{{\bf W}}  \def\mW{{\mathscr W}}
\def\o{\omega}  \def\cX{{\mathcal X}}     \def\bX{{\bf X}}  \def\mX{{\mathscr X}}
\def\x{\xi}     \def\cY{{\mathcal Y}}     \def\bY{{\bf Y}}  \def\mY{{\mathscr Y}}
\def\X{\Xi}     \def\cZ{{\mathcal Z}}     \def\bZ{{\bf Z}}  \def\mZ{{\mathscr Z}}
\def\be{{\bf e}}
\def\bv{{\bf v}} \def\bu{{\bf u}}
\def\Om{\Omega}
\def\bbD{\pmb \Delta}
\def\mm{\mathrm m}
\def\mn{\mathrm n}

\newcommand{\mc}{\mathscr {c}}

\newcommand{\gA}{\mathfrak{A}}          \newcommand{\ga}{\mathfrak{a}}
\newcommand{\gB}{\mathfrak{B}}          \newcommand{\gb}{\mathfrak{b}}
\newcommand{\gC}{\mathfrak{C}}          \newcommand{\gc}{\mathfrak{c}}
\newcommand{\gD}{\mathfrak{D}}          \newcommand{\gd}{\mathfrak{d}}
\newcommand{\gE}{\mathfrak{E}}
\newcommand{\gF}{\mathfrak{F}}           \newcommand{\gf}{\mathfrak{f}}
\newcommand{\gG}{\mathfrak{G}}           
\newcommand{\gH}{\mathfrak{H}}           \newcommand{\gh}{\mathfrak{h}}
\newcommand{\gI}{\mathfrak{I}}           \newcommand{\gi}{\mathfrak{i}}
\newcommand{\gJ}{\mathfrak{J}}           \newcommand{\gj}{\mathfrak{j}}
\newcommand{\gK}{\mathfrak{K}}            \newcommand{\gk}{\mathfrak{k}}
\newcommand{\gL}{\mathfrak{L}}            \newcommand{\gl}{\mathfrak{l}}
\newcommand{\gM}{\mathfrak{M}}            \newcommand{\gm}{\mathfrak{m}}
\newcommand{\gN}{\mathfrak{N}}            \newcommand{\gn}{\mathfrak{n}}
\newcommand{\gO}{\mathfrak{O}}
\newcommand{\gP}{\mathfrak{P}}             \newcommand{\gp}{\mathfrak{p}}
\newcommand{\gQ}{\mathfrak{Q}}             \newcommand{\gq}{\mathfrak{q}}
\newcommand{\gR}{\mathfrak{R}}             \newcommand{\gr}{\mathfrak{r}}
\newcommand{\gS}{\mathfrak{S}}              \newcommand{\gs}{\mathfrak{s}}
\newcommand{\gT}{\mathfrak{T}}             \newcommand{\gt}{\mathfrak{t}}
\newcommand{\gU}{\mathfrak{U}}             \newcommand{\gu}{\mathfrak{u}}
\newcommand{\gV}{\mathfrak{V}}             \newcommand{\gv}{\mathfrak{v}}
\newcommand{\gW}{\mathfrak{W}}             \newcommand{\gw}{\mathfrak{w}}
\newcommand{\gX}{\mathfrak{X}}               \newcommand{\gx}{\mathfrak{x}}
\newcommand{\gY}{\mathfrak{Y}}              \newcommand{\gy}{\mathfrak{y}}
\newcommand{\gZ}{\mathfrak{Z}}             \newcommand{\gz}{\mathfrak{z}}

\def\ve{\varepsilon}   \def\vt{\vartheta}    \def\vp{\varphi}    \def\vk{\varkappa}

\def\A{{\mathbb A}} \def\B{{\mathbb B}} \def\C{{\mathbb C}}
\def\dD{{\mathbb D}} \def\E{{\mathbb E}} \def\dF{{\mathbb F}} \def\dG{{\mathbb G}} \def\H{{\mathbb H}}\def\I{{\mathbb I}} \def\J{{\mathbb J}} \def\K{{\mathbb K}} \def\dL{{\mathbb L}}\def\M{{\mathbb M}} \def\N{{\mathbb N}} \def\O{{\mathbb O}} \def\dP{{\mathbb P}} \def\R{{\mathbb R}}\def\S{{\mathbb S}} \def\T{{\mathbb T}} \def\U{{\mathbb U}} \def\V{{\mathbb V}}\def\W{{\mathbb W}} \def\X{{\mathbb X}} \def\Y{{\mathbb Y}} \def\Z{{\mathbb Z}}


\def\la{\leftarrow}              \def\ra{\rightarrow}            \def\Ra{\Rightarrow}
\def\ua{\uparrow}                \def\da{\downarrow}
\def\lra{\leftrightarrow}        \def\Lra{\Leftrightarrow}


\def\lt{\biggl}                  \def\rt{\biggr}
\def\ol{\overline}               \def\wt{\widetilde}
\def\no{\noindent}


\let\ge\geqslant                 \let\le\leqslant
\def\lan{\langle}                \def\ran{\rangle}
\def\/{\over}                    \def\iy{\infty}
\def\sm{\setminus}               \def\es{\emptyset}
\def\ss{\subset}                 \def\ts{\times}
\def\pa{\partial}                \def\os{\oplus}
\def\om{\ominus}                 \def\ev{\equiv}
\def\iint{\int\!\!\!\int}        \def\iintt{\mathop{\int\!\!\int\!\!\dots\!\!\int}\limits}
\def\el2{\ell^{\,2}}             \def\1{1\!\!1}
\def\sh{\sharp}
\def\wh{\widehat}
\def\bs{\backslash}
\def\intl{\int\limits}

\def\na{\mathop{\mathrm{\nabla}}\nolimits}
\def\sh{\mathop{\mathrm{sh}}\nolimits}
\def\ch{\mathop{\mathrm{ch}}\nolimits}
\def\where{\mathop{\mathrm{where}}\nolimits}
\def\all{\mathop{\mathrm{all}}\nolimits}
\def\as{\mathop{\mathrm{as}}\nolimits}
\def\Area{\mathop{\mathrm{Area}}\nolimits}
\def\arg{\mathop{\mathrm{arg}}\nolimits}
\def\const{\mathop{\mathrm{const}}\nolimits}
\def\det{\mathop{\mathrm{det}}\nolimits}
\def\diag{\mathop{\mathrm{diag}}\nolimits}
\def\diam{\mathop{\mathrm{diam}}\nolimits}
\def\dim{\mathop{\mathrm{dim}}\nolimits}
\def\dist{\mathop{\mathrm{dist}}\nolimits}
\def\Im{\mathop{\mathrm{Im}}\nolimits}
\def\Iso{\mathop{\mathrm{Iso}}\nolimits}
\def\Ker{\mathop{\mathrm{Ker}}\nolimits}
\def\Lip{\mathop{\mathrm{Lip}}\nolimits}
\def\rank{\mathop{\mathrm{rank}}\limits}
\def\Ran{\mathop{\mathrm{Ran}}\nolimits}
\def\Re{\mathop{\mathrm{Re}}\nolimits}
\def\Res{\mathop{\mathrm{Res}}\nolimits}
\def\res{\mathop{\mathrm{res}}\limits}
\def\sign{\mathop{\mathrm{sign}}\nolimits}
\def\span{\mathop{\mathrm{span}}\nolimits}
\def\supp{\mathop{\mathrm{supp}}\nolimits}
\def\Tr{\mathop{\mathrm{Tr}}\nolimits}
\def\BBox{\hspace{1mm}\vrule height6pt width5.5pt depth0pt \hspace{6pt}}


\newcommand\nh[2]{\widehat{#1}\vphantom{#1}^{(#2)}}
\def\dia{\diamond}

\def\Oplus{\bigoplus\nolimits}



\def\qqq{\qquad}
\def\qq{\quad}
\let\ge\geqslant
\let\le\leqslant
\let\geq\geqslant
\let\leq\leqslant
\newcommand{\ca}{\begin{cases}}
\newcommand{\ac}{\end{cases}}
\newcommand{\ma}{\begin{pmatrix}}
\newcommand{\am}{\end{pmatrix}}
\renewcommand{\[}{\begin{equation}}
\renewcommand{\]}{\end{equation}}
\def\eq{\begin{equation}}
\def\qe{\end{equation}}
\def\[{\begin{equation}}
\def\bu{\bullet}

\title[{Spectral band localization for Schr\"odinger operators on
 periodic graphs}]
{Spectral band localization  for Schr\"odinger operators on discrete
periodic graphs}

\date{\today}
\author[Evgeny Korotyaev]{Evgeny Korotyaev}
\address{Mathematical Physics Department, Faculty of Physics, Ulianovskaya 2,
St. Petersburg State University, St. Petersburg, 198904, Russia,
 \ korotyaev@gmail.com,}
\author[Natalia Saburova]{Natalia Saburova}
\address{Department of Mathematical Analysis, Algebra and Geometry, Institute of Mathematics,
Information and Space Technologies, Uritskogo St. 68, Northern (Arctic)
Federal University,
Arkhangelsk, 163002,
 \ n.saburova@gmail.com}


\begin{abstract}
We consider Schr\"odinger operators on periodic discrete  graphs. It
is known that the spectrum of these operators has band structure. We
obtain a localization of spectral bands in terms of eigenvalues of
 Dirichlet and Neumann operators on a finite graph,
which is constructed from the fundamental cell of the periodic
graph. The proof is based on the Floquet decomposition of
Schr\"odinger operators and the minimax principle.
\end{abstract}

\maketitle

\vskip 0.25cm

\section {\lb{Sec1}Introduction and main result}
\setcounter{equation}{0}

 Operators on periodic
graphs are of interest due to their applications to problems of
physics and chemistry. They are used to describe and to study
properties of different periodic media, including nanomedia.
We consider Schr\"odinger operators with periodic potentials on
$\Z^d$-periodic discrete graphs, $d\ge 2$. It is known that the
spectrum of Schr\"odinger operators consists  of an absolutely
continuous part and a finite number of flat bands (i.e., eigenvalues
of infinite multiplicity). The absolutely continuous spectrum
consists of a finite number of intervals (spectral bands) separated
by gaps.  Here we have a well-known problem: to estimate the
spectral bands and gaps in terms of graph parameters and potentials.
In the case of the Schr\"odinger operators $-\D+Q$ with a periodic
potential $Q$ in $\R^d$ there are two-sided estimates of potentials
in terms of gap lengths only at $d=1$ in \cite{K98}, \cite{K03}. We
do not know other estimates.  For the case of periodic graphs we
know only two papers about estimates of spectrum and gaps:

(1) Lled\'o and Post \cite{LP08} considered Laplacians on metric graphs. In
this case they determined estimates (the so-called eigenvalue bracketing)  using
various types of boundary conditions at the vertices. Via an
explicit Cattaneo correspondence \cite{C97} of the equilateral
metric and discrete graph spectra they carry over these estimates
from the metric graph Laplacian to the normalized Laplacian on the
discrete graph. Finally, they wrote {\it "It is a priori not clear
how the eigenvalue bracketing can be seen directly for discrete
Laplacians, so our analysis may serve as an example of how to use
metric graphs to obtain results for discrete graphs"} (p.809 in
\cite{LP08}).

(2) Korotyaev and Saburova \cite{KS13} considered Schr\"odinger
operators on the discrete graphs  and estimated the Lebesgue measure
of their spectrum  in terms of geometric parameters of the graph
only.

\medskip

In our paper we estimate the position of bands of discrete
Schr\"odinger operators on periodic graphs. Here even for the
Laplacian the Lled\'o-Post result \cite{LP08} does not work, since our
Laplacian  is not normalized and the Cattaneo correspondence between
the spectra of Laplacians on discrete and metric graphs treated in
\cite{LP08} does not hold true. We estimate directly spectral band
positions of discrete Schr\"odinger operators in terms of
eigenvalues of Dirichlet and Neumann operators on a finite graph,
which is constructed from the fundamental cell of the periodic
graph. These estimates in some cases allow to determine the
existence of gaps in the spectrum of Schr\"odinger operators. Note
that even for Laplacians it is new.

\subsection{Schr\"odinger operators on periodic graphs.}
Let $\G=(V,\cE)$ be a connected graph, possibly  having loops  and
multiple edges, where $V$ is the set of its vertices and  $\cE$ is
the set of its unoriented edges. An edge connecting vertices $u$
and $v$ from $V$ will be denoted as the unordered pair $(u,v)_e\in
\cE$ and is said to be \emph{incident} to the vertices. Vertices
$u,v\in V$ will be called \emph{adjacent} and denoted by $u\sim v$,
if $(u,v)_e\in \cE$. For each vertex $v\in V$ we define the degree
${\vk}_v=\deg v$ as the number of all its incident edges from $\cE$
(here a loop is counted twice). Below we consider locally finite
$\Z^d$-periodic graphs $\G$, i.e., graphs satisfying the following
conditions:

1) {\it the number of vertices from $V$ in any bounded domain $\ss\R^d$ is
finite;

2) the degree of each vertex is finite;

3) there exists a basis $a_1,\ldots,a_d$ in $\R^d$ such that $\G$ is
invariant under translations through the vectors $a_1,\ldots,a_d$:
$$
\G+a_s=\G, \qqq  \forall\, s\in\N_d=\{1,\ldots,d\}.
$$
The vectors $a_1,\ldots,a_d$ are called the periods of $\G$.}

From
this definition it follows that a $\Z^d$-periodic graph $\G$ is
invariant under translations through any integer vector (in the
basis $a_1,\ldots,a_d$):
$$
\G+\mm=\G,\qqq \forall\, \mm\in\Z^d.
$$

Let $\ell^2(V)$ be the
Hilbert space of all square summable functions $f:V\to \C$, equipped
with the norm
$$
\|f\|^2_{\ell^2(V)}=\sum_{v\in V}|f(v)|^2<\infty.
$$
We  define the self-adjoint Laplacian  (or the Laplace operator)
$\D$ on $f\in\ell^2(V)$ by
\[
\lb{DLO}
 \big(\D f\big)(v)= \sum\limits_{(v,\,u)_e\in\cE}\big(f(v)-f(u)\big), \qqq
 v\in V.
\]

We recall the basic facts  (see \cite{Me94},
\cite{M92}, \cite{MW89}) for both finite and periodic graphs:

 \emph{the point 0 belongs to the
spectrum $\s(\D)$ and $\s(\D)$ is contained in  $[0,2\vk_+]$, i.e.,
}
\[\lb{mp}
0\in\s(\D)\subset[0,2\vk_+],\qqq
\textrm{where}\qqq
\vk_+=\sup_{ v\in V}\deg v<\infty.
\]

We consider the Schr\"odinger operator $H$ acting on the Hilbert
space $\ell^2(V)$ and given by
\[
\lb{Sh}
H=\D+Q,
\]
\[
\lb{Pot}
\big(Q f\big)(v)=Q(v)f(v),\qqq \forall\, v\in V,
\]
where we assume that the potential $Q$ is real valued and satisfies
$$
Q(v+\wt a_s)=Q(v), \qqq  \forall\, (v,s)\in V\ts\N_d,
$$
for some linearly independent integer vectors  $\wt a_1,\ldots,\wt a_d\in\Z^d$ (in the basis $a_1,\ldots,a_d$).
The vectors $\wt a_1,\ldots,\wt a_d$ are called \emph{the periods of the potential} $Q$.
Since the periods $\wt a_1,\ldots,\wt a_d$ of the potential are also periods of the periodic graph, we may assume that the periods of the potential are the same as the periods of the graph.

\subsection{Spectrum of Schr\"odinger operators.}
We define \emph{the
fundamental graph} $\G_F=(V_F,\cE_F)$ of the periodic graph $\G$ as a graph
on the surface $\R^d/\Z^d$ by
\[
\lb{G0} \G_F=\G/{\Z}^d\ss\R^d/\Z^d.
\]
The fundamental graph $\G_F$ has the vertex set $V_F$ and the set
$\cE_F$ of unoriented edges, which are finite. In the space $\R^d$
we consider a coordinate system with the origin at  some point $O$.
The coordinate axes of this system are directed along the vectors
$a_1,\ldots,a_d$. Below the coordinates of all vertices of $\G$ will
be expressed  in this coordinate system. We identify the vertices of
the fundamental graph $\G_F=(V_F,\cE_F)$ with the vertices of the
graph $\G=(V,\cE)$ from the set $[0,1)^d$ by
\[
\lb{V0} V_F=[0,1)^d\cap V=\{v_1,\ldots,v_\n\},\qqq \n=\# V_F<\infty,
\]
where $\n=\# V_F$ is the number of vertices of $\G_F$.

The Schr\"odinger operator $H=\D+Q$ on $\ell^2(V)$ has  the
decomposition into a constant fiber direct integral
\[
\lb{raz}
\begin{aligned}
& \ell^2(V)={1\/(2\pi)^d}\int^\oplus_{\T^d}\ell^2(V_F)\,d\vt ,\qqq
UH U^{-1}={1\/(2\pi)^d}\int^\oplus_{\T^d}H(\vt)d\vt,
\end{aligned}
\]
$\T^d=\R^d/(2\pi\Z)^d$,     for some unitary operator $U$. Here
$\ell^2(V_F)=\C^\nu$ is the fiber space  and the Floquet $\nu\ts\nu$
matrix  $H(\vt)$ (i.e., a fiber  matrix) is given by
\[
\lb{Hvt}
H(\vt)=\D(\vt)+q,\qqq q=\diag(q_1,\ldots,q_\n),\qqq \forall\,\vt\in \T^d,
\]
and $q_j$ denote the values of the potential $Q$ on the vertex set
$V_F$ by
\[\lb{pott}
Q(v_j)=q_j, \qqq j\in \N_\n=\{1,\ldots,\nu\}.
\]
The decomposition \er{raz} is standard and follows from the
Floquet-Bloch theory \cite{RS78}. The precise expression of the
Floquet matrix $\D(\vt)=\{\D_{jk}(\vt)\}_{j,k=1}^\n$ for the
Laplacian $\D$ is given by \er{l2.15'}. Each Floquet ${\nu\ts\nu}$
matrix $H(\vt)$, $\vt\in\T^d$, has $\n$ eigenvalues labeled by
$$
\l_1(\vt)\leq\ldots\leq\l_\n(\vt).
$$
Note that the spectrum of the Floquet matrix $H(\vt)$ does not
depend  on the choice of the coordinate origin $O$. Each
$\l_n(\cdot)$, $n\in\N_\n$, is a real and continuous function on the
torus $\T^d$ and creates the spectral band $\s_n(H)$ given by
\[
\lb{ban} \s_n=\s_n(H)=[\l_n^-,\l_n^+]=\l_n(\T^d), \qqq
\l_n^-=\min_{\vt \in\T^d}\l_n(\vt),\qqq \l_n^+=\max_{\vt
\in\T^d}\l_n(\vt ).
\]
Thus, the spectrum of the operator $H$ on the periodic graph $\G$ is
given by
\[
\lb{r0}
\s(H)=\bigcup_{\vt\in\T^d}\s\big(H(\vt)\big)=\bigcup_{n=1}^{\nu}\s_n(H).
\]
Note that if $\l_n(\cdot)= C_n=\const$ on some set $\mB\ss\T^d$ of
positive Lebesgue measure, then  the operator $H$ on $\G$ has the
eigenvalue $C_n$ with infinite multiplicity. We call $C_n$ a
\emph{flat band}. Thus, the spectrum of the Schr\"odinger operator
$H$ on the periodic graph $\G$ has the form
\[
\lb{r0}
\s(H)=\s_{ac}(H)\cup \s_{fb}(H),
\]
where $\s_{ac}(H)$ is the absolutely continuous spectrum, which is a
union  of non-degenerated intervals, and $\s_{fb}(H)$ is the set of
all flat bands (eigenvalues of infinite multiplicity). An open
interval between two neighboring non-degenerated spectral bands is
called a \emph{spectral gap}.

\subsection{Localization of spectral bands for Schr\"odinger operators.}
In order to formulate estimates of spectral band positions we define
two operators on a finite graph $\G_N=(V_N,\cE_N)$ with a vertex set
$V_N$ and an edge set $\cE_N$. This graph is introduced by the
following way. Denote by $\cE_N^{in}$ the set of all edges of $\G$
connecting the vertices from $V_F$, defined by \er{V0}. An edge
connecting a vertex from $V_F$ with a vertex from $V\setminus V_F$
is called a \emph{bridge}.
Denote by $\cB$ the set of all bridges. Since $\Z^d$-periodic graph $\G$ is connected, the number $\b=\#\cB$ of all bridges of $\G$ depends on the choice of the coordinate origin $O$ and satisfies $2d\leq\b$. Due to periodicity of the graph the set $\cB$ consists of pairs of
equivalent to each other (with respect to the action of the group
$\Z^d$) bridges. We take one bridge from each pair and denote the obtained set of bridges by $\cB_N$. The finite graph $\G_N=(V_N,\cE_N)$ is the
edge-induced subgraph  of $\G$,  its edge set is
$\cE_N=\cE_N^{in}\cup\cB_N$ and its vertex set $V_N$ consists of all
ends of edges of $\cE_N$, i.e.,
\[
\lb{VN} \cE_N=\cE_N^{in}\cup\cB_N,\qqq V_N=V_F\cup \{u\in V:
(u,v)_e\in \cB_N, v\in V_F \}.
\]
Let $\vk_v^N=\deg v$ be the degree of the vertex $v\in V_N$ on the
graph $\G_N$. A vertex $v\in V_N$ will be called an \emph{inner}
vertex of $\G_N$, if $\vk_v=\vk_v^N$, i.e., if all its adjacent
(neighbor) vertices from $V$ also belong to the graph $\G_N$. Denote
by $V_{D}$ the set of all inner vertices of $\G_N$. Let
$\n_\f=\# V_\f$  be the number of the vertices in $V_\f, \ \f =D,N$. We define a
\emph{boundary} $\pa V_N$ of $\G_N$ by the standard identity:
\[\lb{inn}
\pa V_N=V_N\sm V_{D}.
\]
\begin{figure}[h]
\centering
\unitlength 1.0mm 
\linethickness{0.4pt}
\ifx\plotpoint\undefined\newsavebox{\plotpoint}\fi 

\begin{picture}(120,55)(0,0)
\unitlength 1.1mm
\put(41,48){$\scriptstyle a_2$}
\put(62.5,27.5){$\scriptstyle a_1$}
\put(42.5,28){$\scriptstyle O$}
\put(20,8){\emph{a)}}

\put(56,19.5){$\scriptstyle v_6$}
\put(56,40){$\scriptstyle v_3$}
\put(57.0,28.3){$\scriptstyle v_1$}

\put(66.5,19.5){$\scriptstyle v_7$}
\put(41.0,20.0){$\scriptstyle v_5$}
\put(66.5,39.5){$\scriptstyle v_4$}
\put(41.0,40.0){$\scriptstyle v_2$}

\put(45.0,30){\vector(1,0){20.00}}
\put(45.0,30){\vector(0,1){20.00}}
\multiput(45,50)(4,0){5}{\line(1,0){2}}
\multiput(65,30)(0,4){5}{\line(0,1){2}}

\put(35,10){\line(0,1){40.00}}
\put(55,10){\line(0,1){40.00}}
\put(75,10){\line(0,1){40.00}}

\put(35,10){\line(-1,1){10.00}}
\put(35,10){\line(1,1){10.00}}
\put(35,30){\line(1,-1){10.00}}
\put(35,30){\line(-1,-1){10.00}}

\put(35,30){\line(-1,1){10.00}}
\put(35,30){\line(1,1){10.00}}
\put(35,50){\line(1,-1){10.00}}
\put(35,50){\line(-1,-1){10.00}}

\put(75,10){\line(-1,1){10.00}}
\put(75,10){\line(1,1){10.00}}
\put(75,30){\line(1,-1){10.00}}
\put(75,30){\line(-1,-1){10.00}}

\put(75,30){\line(-1,1){10.00}}
\put(75,30){\line(1,1){10.00}}
\put(75,50){\line(1,-1){10.00}}
\put(75,50){\line(-1,-1){10.00}}

\put(55,10){\line(-1,1){10.00}}
\put(55,10){\line(1,1){10.00}}
\put(55,30){\line(1,-1){10.00}}
\put(55,30){\line(-1,-1){10.00}}

\put(55,30){\line(-1,1){10.00}}
\put(55,30.2){\line(-1,1){10.00}}
\put(55,29.8){\line(-1,1){10.00}}

\put(55,30){\line(1,1){10.00}}
\put(55,30.2){\line(1,1){10.00}}
\put(55,29.8){\line(1,1){10.00}}

\put(55,30){\line(1,-1){10.00}}
\put(55,30.2){\line(1,-1){10.00}}
\put(55,29.8){\line(1,-1){10.00}}

\put(55,50.0){\line(1,-1){10.00}}
\put(55,50){\line(-1,-1){10.00}}

\put(55,30.2){\line(-1,-1){10.00}}
\put(55,30){\line(-1,-1){10.00}}
\put(55,29.8){\line(-1,-1){10.00}}

\put(35,10){\circle{1}}
\put(55,10){\circle{1}}
\put(75,10){\circle{1}}

\put(35,30){\circle{1}}
\put(55,30){\circle{1}}
\put(75,30){\circle{1}}

\put(35,50){\circle{1}}
\put(55,50){\circle{1}}
\put(75,50){\circle{1}}

\put(25,20){\circle{1}}
\put(35,20){\circle{1}}
\put(45,20){\circle*{1}}
\put(55,20){\circle*{1}}
\put(55.1,20){\line(0,1){20.00}}
\put(55.2,20){\line(0,1){20.00}}
\put(55.0,20){\line(0,1){20.00}}
\put(54.9,20){\line(0,1){20.00}}
\put(54.8,20){\line(0,1){20.00}}
\put(65,20){\circle*{1}}
\put(75,20){\circle{1}}
\put(85,20){\circle{1}}

\put(25,40){\circle{1}}
\put(35,40){\circle{1}}
\put(45,40){\circle*{1}}
\put(55,40){\circle*{1}}
\put(65,40){\circle*{1}}
\put(75,40){\circle{1}}
\put(85,40){\circle{1}}

\put(55,30){\circle*{1}}
\put(65,40){\circle*{1}}
\put(55,50){\circle{1}}
\end{picture}

\begin{picture}(130,40)(0,0)
\unitlength 1.0mm
\put(0,7){\vector(1,0){130.00}}
\put(0,7.1){\line(1,0){20}}
\put(0,6.9){\line(1,0){20}}
\put(0,7.2){\line(1,0){20}}
\put(0,6.8){\line(1,0){20}}
\put(20,6){\line(0,1){2.00}}

\put(25,7.1){\line(1,0){15}}
\put(25,6.9){\line(1,0){15}}
\put(25,7.2){\line(1,0){15}}
\put(25,6.8){\line(1,0){15}}
\put(25,6){\line(0,1){2.0}}
\put(40,6){\line(0,1){2.0}}

\put(60,6.9){\line(1,0){35}}
\put(60,6.8){\line(1,0){35}}
\put(60,7.0){\line(1,0){35}}
\put(60,7.2){\line(1,0){35}}
\put(60,7.1){\line(1,0){35}}
\put(60,5.8){\line(0,1){2.0}}
\put(95,5.8){\line(0,1){2.0}}

\put(0,6){\line(0,1){2.00}}
\put(120,6){\line(0,1){2.00}}
\put(-0.5,3){$\scriptstyle 0$}
\put(31.5,4){$\scriptstyle \s_2$}
\put(75,4){$\scriptstyle \s_3$}
\put(119,3){$\scriptstyle 12$}
\put(125.0,3){$\scriptstyle \s(\D)$}

\put(9.0,4){$\scriptstyle \s_1$}

\put(0,18){\vector(1,0){130.00}}
\put(0,17){\line(0,1){2.00}}
\put(120,17){\line(0,1){2.00}}

\put(122.0,20){$\scriptstyle \s(H_D)$}

\put(60,18){\circle*{1}}
\bezier{25}(60,7)(60,15.5)(60,24)

\put(0,24.1){\line(1,0){60}}
\put(0,24.2){\line(1,0){60}}
\put(0,24.3){\line(1,0){60}}
\put(0,24.4){\line(1,0){60}}
\put(0,24.5){\line(1,0){60}}

\put(20,22.5){\line(1,0){100}}
\put(20,22.4){\line(1,0){100}}
\put(20,22.3){\line(1,0){100}}
\put(20,22.2){\line(1,0){100}}
\put(20,22.1){\line(1,0){100}}

\put(25,23.1){\line(1,0){95}}
\put(25.0,23.2){\line(1,0){95}}
\put(25.0,23.3){\line(1,0){95}}
\put(25.0,23.4){\line(1,0){95}}
\put(25.0,23.5){\line(1,0){95}}

\put(60,25.1){\line(1,0){35}}
\put(60,25.2){\line(1,0){35}}
\put(60,25.3){\line(1,0){35}}
\put(60,25.4){\line(1,0){35}}
\put(60,25.5){\line(1,0){35}}
\put(10.0,21.0){$\scriptstyle J_1$}

\put(43.0,19.5){$\scriptstyle J_2$}

\put(10.6,26.5){$\scriptstyle \wt J_1=\wt J_2$}

\put(76.0,26.5){$\scriptstyle\wt J_3$}
\put(105.5,24.3){$\scriptstyle J_3$}

\put(0,25.5){\line(1,0){40.00}}
\put(0,25.4){\line(1,0){40.00}}
\put(0,25.3){\line(1,0){40.00}}
\put(0,25.2){\line(1,0){40.00}}
\put(0,25.1){\line(1,0){40.00}}

\put(0,30){\vector(1,0){130.00}}
\put(120,29){\line(0,1){2.00}}

\put(122.0,32){$\scriptstyle \s(H_N)$}
\put(0,30){\circle*{1}}
\put(20,30){\circle*{1}}
\put(25,30){\circle*{1}}
\put(40.2,31.5){\circle*{1}}
\put(40.2,30){\circle*{1}}
\put(40.2,28.5){\circle*{1}}
\put(95,30){\circle*{1}}

\bezier{25}(0,30)(0,18.5)(0,7)
\bezier{30}(20,30)(20,18.5)(20,7)
\bezier{10}(25,30)(25,27)(25,24)
\bezier{30}(40,30)(40,18.5)(40,7)
\bezier{30}(95,30)(95,18.5)(95,7)

\put(-9,5){\emph{b)}}

\put(0,10.5){\line(1,0){40}}
\put(0,10.4){\line(1,0){40}}
\put(0,10.3){\line(1,0){40}}
\put(0,10.2){\line(1,0){40}}
\put(0,10.1){\line(1,0){40}}
\put(7.0,11.5){$\scriptstyle J_1\cap\wt J_1$}

\put(20,12.1){\line(1,0){20}}
\put(20,12.2){\line(1,0){20}}
\put(20,12.3){\line(1,0){20}}
\put(20,12.4){\line(1,0){20}}
\put(20,12.5){\line(1,0){20}}

\put(27.0,13.5){$\scriptstyle J_2\cap\wt J_2$}
\put(60,10.1){\line(1,0){35}}
\put(60,10.2){\line(1,0){35}}
\put(60,10.3){\line(1,0){35}}
\put(60,10.4){\line(1,0){35}}
\put(60,10.5){\line(1,0){35}}
\put(73.3,11.5){$\scriptstyle J_3\cap\wt J_3$}
\end{picture}
\caption{ \footnotesize \emph{a)} A periodic graph $\G$ and its finite graph
$\G_N$,  the vertices of the graph $\G_N$ are black; the edges of
$\G_N$ are marked by bold lines. The set of the inner vertices and the boundary  are $V_D=\{v_1\}$ and $\partial V_N=\{v_2,v_3,v_4,v_5,v_6,v_7\}$, respectively.
\emph{b)} Eigenvalues of the operators $H_N$ and $H_D$, the intervals $J_n$ and $\wt J_n$, $n\in\N_3$, and their intersections, the spectrum of the Laplacian $\D$.} \label{ff.0.11}
\end{figure}
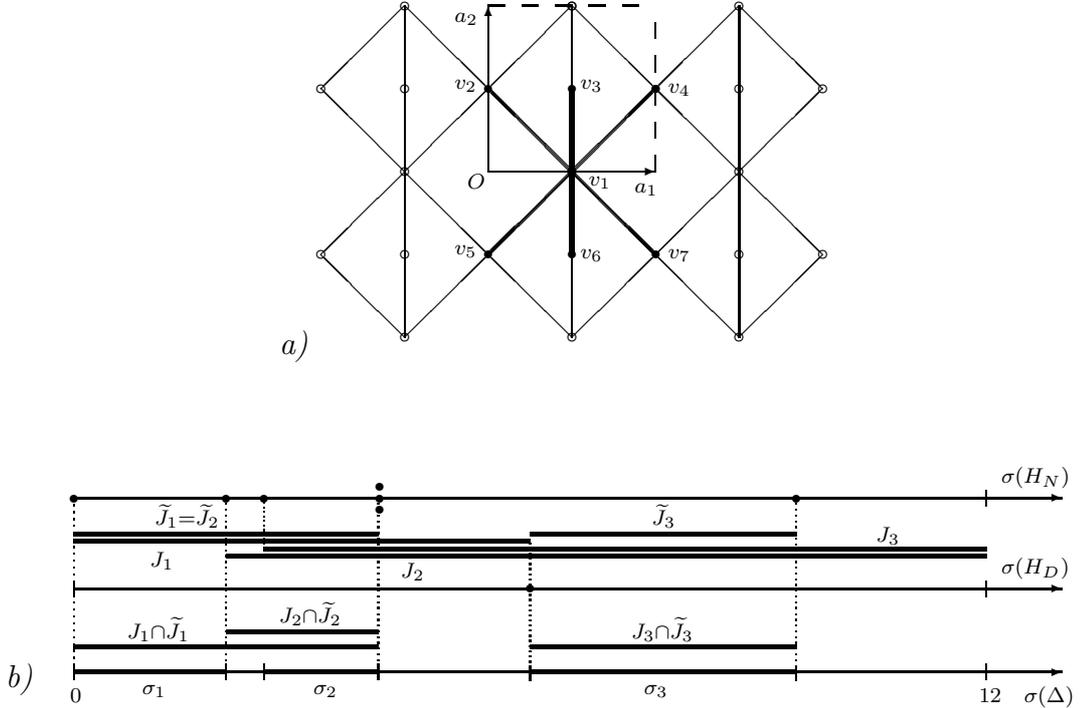

{\bf Example.} For $\Z^2$-periodic
graph shown in Fig.\ref{ff.0.11}\emph{a} the set $V_F$ consists of the
vertices \linebreak $\{v_1,v_2,v_3\}$. The finite graph
$\G_N=(V_N,\cE_N)$ has the vertex set $V_N$ given by
$$
V_N=\{v_1,v_2,v_3,v_4=v_2+a_1,v_5=v_2-a_2,v_6=v_3-a_2,v_7=v_2+a_1-a_2\}.
$$
The set of the inner vertices $V_D$ and the boundary
$\partial V_N$ of $\G_N$ have the form
$$
V_D=\{v_1\},\qqq \partial V_N=\{v_2,v_3,v_4,v_5,v_6,v_7\}.
\qqq \BBox
$$

On the finite graph $\G_N$ we define two self-adjoint operators
$H_N$ and $H_D$:

1) The Neumann operator $H_N$ on $\ell^2(V_N)$ is defined by
\[
\lb{Sh+}
H_N=\D_N+Q,
\]
where
\[
\lb{DLO+} \big(\D_N f\big)(v)=\r_v\vk_v^Nf(v)-
\sqrt{\r_v}\sum\limits_{(v,\,u)_e\in\cE_N}\sqrt{\r_u}\,f(u), \qqq
 v\in V_N, \qq f\in \ell^2(V_N),
\]
\[
\lb{Vj+} \r_v=\# V_v,\qqq  V_v=\big(\{v\}+\Z^d\big)\cap V_N, \qqq
v\in V_N,
\]
i.e., $\r_v$ is the number of vertices from $V_N$ equivalent to $v$
with respect to the action of the group $\Z^d$, $1\leq\r_v\leq\n_N$.
Note that $\vk_j^N=\vk_j$ and $\r_j=1$ for all $j\in\N_{\n_D}$.

2) The Dirichlet operator $H_D$ on $f\in \ell^2(V_N)$ is also
defined by \er{Sh+}, but with the Dirichlet boundary conditions
$f(v)=0$ for all $v\in\pa V_N$. We will identify the Dirichlet
operator $H_D$ on $f\in \ell^2(V_N)$ with the Dirichlet operator
$H_D$ on $f\in \ell^2(V_D)$, since $f(v)=0$ for all $v\in\pa V_N$.

Denote the eigenvalues of the operators $H_\f$, $\f=D,N$, counted
according to multiplicity, by
\[
\lb{LDN} \l_1^\f \le \l_2^\f \le\ldots\le \l_{\n_\f}^\f,\qqq
\n_\f=\# V_\f,\qq \f=D,N.
\]
 We rewrite the sequence
$q_1,\ldots,q_\n$ define by \er{pott} in nondecreasing order
\[
\lb{wtqn}
q^\bullet_1\le q^\bullet_2 \le\ldots \le  q^\bullet_\n.
 \]
Here
$q^\bullet_1=q_{n_1},q^\bullet_2=q_{n_2},\ldots,q^\bullet_\n=q_{n_\n}$
for some distinct numbers $n_1, {n_2},\ldots,{n_\n}\in\N_\n$.

\begin{theorem}
\label{T1} Each spectral band $\s_n(H)$ of the operator $H=\D+Q$ on
the graph $\G$ satisfies
\[\lb{ell}
\s_n(H)\subset J_n\cap \wt J_n,\qqq n\in\N_\n,
\]
where the intervals $J_n$, $\wt J_n$ are given by
\[
\lb{J1}
J_n=\ca [\l_n^N,\l_n^D],  & n=1,\ldots, \n_D,\\[2pt]
[\l_n^N,q^\bu_n+2\vk_+],  & n=\n_D+1,\ldots,\n, \ac
\]
and
\[\lb{wJ1}
\wt J_n=\ca [q^\bullet_n,\l_{n+\n_N-\n}^N],  & n=1,\ldots,\n-\n_D,\\[2pt]
[\l_{n-\n+\n_D}^D,\l_{n+\n_N-\n}^N], & n=\n-\n_D+1,\ldots,\n \ac.
\]
\end{theorem}

\no \textbf{Remark.} 1) Due to  Cattaneo correspondence \cite{C97}
Lled\'o and Post \cite{LP08} considered the so-called normalized
Laplacian $\hat\D$ on $\ell^2(V)$ given by \ $ \hat\D=\c \D \c $, \
where $\c$ is the multiplication operator on $\ell^2(V)$ given by
$(\c f)(v)=\vk_v^{-{1\/2}}f(v)$ and $\vk_v^{{1\/2}}>0$.
 They estimated the position of the band $\s_n(\hat\D)$ for the normalized
Laplacian $\hat\D$ by
$ \s_n(\hat\D)\subset \hat J_n,\  n\in\N_\n, $ where the segments
$\hat J_n$ have the form similar to $J_n$ from \er{J1}.

2) Let the graph $\G$ be bipartite and regular of degree $\vk_+$, i.e., each its vertex $v$ has
the degree $\vk_v=\vk_+$. If
$H=\D$, then $\wt J_n=\zeta(J_n)$  for each $n\in\N_\n$, where
$\zeta(z)=2\vk_+-z$. Thus, in this case the estimate \er{ell} has
the form
$$
\s_n(\D)\subset J_n\cap\zeta(J_n),\qqq n\in\N_\n.
$$

3) Theorem \ref{T1} estimates the positions of the spectral bands in
terms  of  eigenvalues of the operators $H_N$ and $H_D$ on the
finite graph $\G_N$. Moreover, in some cases it allows to detect the
existence of gaps in the spectrum of the Schr\"odinger operator $H$.
For example, for the graph shown in Fig.\ref{ff.0.11}\emph{a} in the
case when $H=\D$ the intervals $J_n\cap \wt J_n$, $n\in\N_\n$, are
shown in Fig.\ref{ff.0.11}\emph{b}. The spectrum of the Laplacian
$\D$ is also shown in this figure. As we can see Theorem~\ref{T1}
detects \emph{precisely} the second gap in the spectrum of the operator (for more
details see Subsection \ref{S3}).

Now we estimate the total length of all spectral bands of
$H$.

\begin{theorem}
\label{T2} i) The total length of all spectral bands $\s_n(H)$, $n\in\N_\n$, of $H$
satisfies
\[\lb{est1}
\begin{aligned}
\sum_{n=1}^{\n}|\s_n(H)|
\le\sum_{n=\n_D+1}^{\n}(q_n^\bu+2\vk_+-h_n)+
\sum_{n=\n+1}^{\n_N}\l_n^N,
\end{aligned}
\]
where $h_n=\r_n(\vk_n-\vk_{nn}+q_{n})$, $\r_n=\# V_n$,
$\vk_{nn}$ is the number of loops in the vertex $v_n$ on the graph $\G$;
\[\lb{est2}
\sum_{n=1}^{\n}|\s_n(H)|
\le\sum_{n=1}^{\n-\n_D}\big(\l_{\n_N-(\n-\n_D)+n}^N-\l_n^N\big).
\]

ii) The numbers $\n_\f=\# V_\f, \ \f =D,N$ satisfy
\[\lb{nv}
\textstyle 0\leq\n_D\leq\n-1,\qqq \n+d\leq\n_N\leq\n+{\b\/2}\,,
\]
where $\b=\#\cB$ is the number of the bridges of $\G$. Moreover, the boundaries of the inequalities \er{nv} are achieved.
\end{theorem}

\no \textbf{Remark.} For the global estimate of the Lebesgue measure of the spectrum of Schr\"odinger operators $H$ it is enough to know the eigenvalues of the Neumann operator $H_N$.

\section{Proof of the main result}
\setcounter{equation}{0}

\subsection{The Floquet matrix for the Schr\"odinger operator.}
We need to introduce the two oriented edges $(u,v)$ and $(v,u)$ for
each unoriented edge $(u,v)_e\in \cE$: the oriented edge starting at
$u\in V$ and ending at $v\in V$ will be denoted as the ordered pair
$(u,v)$. We denote the sets of all oriented edges of the graph $\G$
and the fundamental graph $\G_F$ by $\cA$ and $\cA_F$, respectively.

We introduce {\it an edge index}, which is important to study the
spectrum of Schr\"odinger operators on periodic graphs. For any
$v\in V$ the following unique representation holds true:
\[
\lb{Dv} v=[v]+\tilde v, \qquad [v]\in\Z^d,\qquad \tilde v\in
V_F\subset[0,1)^d.
\]
In other words, each vertex $v$ can be represented uniquely as the
sum  of an integer part $[v]\in \Z^d$ and a fractional part $\tilde
v$ that is a vertex of $V_F$ defined in \er{V0}. For any oriented
edge $\be=(u,v)\in\cA$ we define {\bf the edge "index"}  $\t({\bf
e})$ as the integer vector
\[
\lb{in} \t({\bf e})=[v]-[u]\in\Z^d,
\]
where due to \er{Dv} we have
$$
u=[u]+\tilde{u},\qquad v=[v]+\tilde{v}, \qquad [u],
[v]\in\Z^d,\qquad \tilde{u},\tilde{v}\in V_F.
$$

If $\be=(u,v)$ is an oriented edge of the graph $\G$, then by the
definition of the fundamental graph there is an oriented edge
$\tilde\be=(\tilde u,\tilde v\,)$ on $\G_F$. For the edge
$\tilde\be\in\cA_F$ we define the edge index $\t(\tilde{\bf e})$ by
\[
\lb{inf} \t(\tilde{\bf e})=\t(\be).
\]
In other words, edge indices of the fundamental graph $\G_F$  are
induced by edge indices of the periodic graph $\G$. The edge
indices, generally speaking, depend on the choice of the coordinate
origin $O$. But in a fixed coordinate system the index of the
fundamental graph edge is uniquely determined by \er{inf}, since
$$
\t(\be+\mm)=\t(\be),\qqq \forall\, (\be,\mm)\in\cA \ts \Z^d.
$$
Note that all bridges of the
graph $\G$ have nonzero indices.

The Schr\"odinger operator $H=\D+Q$ acting on $\ell^2(V)$ has the decomposition
into a constant fiber direct integral \er{raz},
where the Floquet $\nu\ts\nu$
matrix  $H(\vt)$ has the form
\[
\lb{Hvt}
H(\vt)=\D(\vt)+q,\qqq q=\diag(q_1,\ldots,q_\n),\qqq \forall\,\vt\in \T^d.
\]
The Floquet matrix
$\D(\vt)=\{\D_{jk}(\vt)\}_{j,k=1}^\n$ for the Laplacian $\D$ is given by
\medskip
\[
\label{l2.15'}
\D_{jk}(\vt )=\vk_j\d_{jk}-\ca
\sum\limits_{{\bf e}=(v_j,\,v_k)\in{\cA}_F}e^{\,i\lan\t
({\bf e}),\,\vt\ran }, \qq &  {\rm if}\  \ (v_j,v_k)\in \cA_F \\
\qqq 0, &  {\rm if}\  \ (v_j,v_k)\notin \cA_F
\ac,
\]
see \cite{KS13}, where $\vk_j$ is the degree of $v_j$, $\d_{jk}$ is
the Kronecker delta and $\lan\cdot\,,\cdot\ran$ denotes the standard
inner product in $\R^d$.

Now we need the following simple fact (see Theorem 4.3.1 in
\cite{HJ85}). \emph{Let $A,B$ be self-adjoint $\nu\ts\nu$ matrices.
Denote by $\l_1(A)\leq\ldots\leq\l_\n(A)$,
$\l_1(B)\leq\ldots\leq\l_\n(B)$ the eigenvalues of  $A$ and $B$,
respectively, arranged in increasing order, counting multiplicities.
Then we have}
\[
\lb{MP} \l_n(A)+\l_1(B)\leq\l_n(A+B)\leq\l_n(A)+\l_\n(B) \qqq
\forall\, n\in\N_\n.
\]

Inequalities \er{MP} and the basic fact \er{mp} give that the eigenvalues of the Floquet matrix $H(\vt)$ for the
Schr\"odinger operator $H=\D+Q$, satisfy
\[
\lb{snH}
\begin{aligned}
q_n^\bullet\le\l_n(\vt)\le q_n^\bullet+2\vk_+,\qqq \forall\, (\vt,n)\in\T^d\ts\N_\n,
\\
\s_n(H)=\l_n(\T^d)\subset[q_n^\bullet,q_n^\bullet+2\vk_+], \qqq \forall\, n\in\N_\n.
\end{aligned}
\]

\subsection{Proof of the main result.}
Without loss of generality we may assume that the set of inner vertices $V_D$ of the graph ${\G_N=(V_N,\cE_N)}$ has the form
$$
V_D=\{v_1,\ldots,v_{\n_D}\}.
$$
We denote the equivalence classes from $V_N/\Z^d$ by
\[\lb{Vj}
V_j=\big(\{v_j\}+\Z^d\big)\cap V_N, \qqq j\in\N_\n.
\]

The Neumann operator $H_N$ on the graph $\G_N$ is equivalent to the $\n_N\ts\n_N$ self-adjoint matrix
$H_N=\{H_{jk}^N\}_{j,k=1}^{\n_N}$ given by
\[
\lb{H+}
H_N=\D_N+q_N,\qqq q_N=\diag(q_1^N,\ldots,q_{\n_N}^N),
\]
where $q^N_k=q_j$, if $v_k\in V_j$, $k\in\N_{\n_N}$, $j\in\N_\n$,
and the matrix
$\D_N=\{\D_{jk}^N\}_{j,k=1}^{\n_N}$ has the form
\[
\label{l2.15}
\D_{jk}^N=\sqrt{\r_j\r_k}\,(\vk_j^N\d_{jk}-\vk^N_{jk}).
\]
Here $\vk^N_j$ is the degree of the vertex $v_j\in V_N$ on the graph
$\G_N$, $\vk^N_{jk}\ge 1$ is the multiplicity of the edge
$(v_j,v_k)\in\cE_N$ and  $\vk_{jk}^N=0$ if $(v_j,v_k)\notin\cE_N$,
$\r_j=|V_j|$ is the number of vertices in $V_j$.

The Dirichlet operator $H_D$ is described by the $\n_D\ts\n_D$
self-adjoint matrix $H_D=\{H_{jk}^D\}_{j,k=1}^{\n_D}$ with entries
\[\lb{H-}
H_{jk}^D=H_{jk}^N\qqq
\textrm{for all}\qqq j,k\in\N_{\n_D}.
\]

Recall that 
\[\lb{rc}
\vk_j^N=\vk_j \ \textrm{ and } \ \r_j=1 \ \textrm{ for all } \ j\in\N_{\n_D}.
\]

\no {\bf Proof of Theorem \ref{T1}.} For each $\vt\in\T^d$ we define the $\n$-dimensional
subspace $Y_\vt$ of $\C^{\n_N}$ by
\begin{multline}\lb{Yt}
Y_\vt=\big\{x=(x_k)_{k=1}^{\n_N}\in\C^{\n_N} : \forall \,
k=\n+1,\ldots,\n_N \qq x_k=\textstyle e^{\,i\lan
v_k-v_j,\,\vt\ran}\,{x_j}\,, \\ \textrm{where $j=j(k)\in\N_\n$ is such that $v_k\in V_j$}\big\}.
\end{multline}
Note that $j=j(k)$ in \er{Yt} is uniquely defined for each $k=\n+1,\ldots,\n_N$.

Let $X$ be the
$\n_D$-dimensional subspace of $\C^{\n}$, defined by
$$
X=\{x\in\C^{\n} : x_{\n_D+1}=\ldots=x_\n=0\}.
$$

We recall well-known facts.

\emph{Denote by $\l_1(A)\leq\ldots\leq\l_\n(A)$ the eigenvalues of a
self-adjoint $\n\ts\n$ matrix $A$, arranged in increasing order,
counting multiplicities. Each $\l_n$ satisfies the minimax principle:}
\[
\lb{CF1} \l_n(A)=\min_{S_n\subset\C^\n}\max_{\|x\|=1 \atop x\in
S_n}\lan Ax,x\ran,
\]
\[
\lb{CF2} \l_n(A)=\max_{S_{\n-n+1}\subset\C^\n}\min_{\|x\|=1 \atop
x\in S_{\n-n+1}}\lan Ax,x\ran,
\]
\emph{where $S_n$ denotes a subspace of dimension $n$ and the outer
optimization is over all subspaces of the indicated dimension} (see
p.180 in \cite{HJ85}).

First, let $1\leq n\leq\n$. Using \er{CF1} and \er{CF2} we write
\[
\lb{CF3'} \l_j^N=\min_{S_j\ss \C^{\n_N}} \max_{\|x\|=1 \atop
{x\in S_j}}\lan H_Nx,x\ran\geq\min_{S_j\ss \C^{\n_N}} \max_{\|x\|=1
\atop {x\in S_j\cap Y_\vt}}\lan H_Nx,x\ran, \qq j=n+\n_N-\n,
\]
\[
\lb{CF3} \l_n^N=\max_{S_k\ss\C^{\n_N}}\min_{\|x\|=1 \atop {x\in
S_k}}\lan H_Nx,x\ran\leq\max_{S_k\ss\C^{\n_N}}\min_{\|x\|=1 \atop
x\in S_k\cap Y_\vt}\lan H_Nx,x\ran, \qq k=\n_N-n+1,
\]
where $S_j$ denotes a subspace of dimension $j$. For $x\in Y_\vt$ we have
\[\lb{CF5.11}
\lan H_Nx,x\ran=\sum_{j,k=1}^{\n_N}H^N_{jk}\,\bar
x_j\,x_k=\sum_{j=1}^{\n_N}(\r_j\vk_j^N+q_j^N)|x_j|^2-
\sum_{j,k=1}^{\n_N}\vk_{jk}^N\,\sqrt{\r_j\r_k}\,\bar
x_j\,x_k,
\]
where
\begin{multline}\lb{CF5.22}
\sum_{j=1}^{\n_N}(\r_j\vk_j^N+q^N_j)|x_j|^2=
\sum_{j=1}^{\n_D}(\vk_j+q_j)|x_j|^2+\sum_{j=\n_D+1}^{\n}|x_j|^2
\sum\limits_{v\in
V_j}(\r_j\vk_v^N+q_j)\\=\sum_{j=1}^{\n_D}(\vk_j+q_j)|x_j|^2
+\sum_{j=\n_D+1}^{\n}\r_j(\vk_j+q_j)|x_j|^2,
\end{multline}
\[\lb{CF5.33}
\sum_{j,k=1}^{\n_N}\vk_{jk}^N\,\sqrt{\r_j\r_k}\;\bar
x_j\,x_k=
\sum_{j,k=1}^{\n}\sqrt{\r_j\r_k}\sum\limits_{{\bf
e}=(v_j,\,v_k)\in\cA_F}e^{\,i\lan\t ({\bf e}),\,\vt\ran }\,\bar
x_j\,x_k.
\]
In \er{CF5.22} we have used the identities \er{rc} and
\[\lb{deg}
\sum\limits_{v\in V_j}\vk_v^N=\vk_j.
\]

We introduce the new vector
\[\lb{vy}
y=(y_j)_{j=1}^\n,\qqq y_j=\sqrt{\r_j}\;x_j, \qqq j\in\N_\n.
\]
Since $\vk_j=1$ for $1\leq j\leq \n_D$, we have $y_j=x_j$, $j\in\N_{\n_D}$, and, using \er{deg}, for
$x\in Y_\vt$ we have
\[\lb{CF6}
\|x\|^2=\sum_{j=1}^{\n_N}|x_j|^2=\sum_{j=1}^{\n_D}|x_j|^2+\sum_{j=\n_D+1}^{\n}
\r_j\,|x_j|^2=\sum_{j=1}^{\n}|y_j|^2=\|y\|^2.
\]
Combining \er{CF5.11} -- \er{CF5.33} for $x\in Y_\vt$, \er{vy} and the
definition of  $H(\vt)$ in \er{Hvt} we obtain
\[\lb{CF5.111}
\lan H_Nx,x\ran=\sum_{j=1}^{\n}(\vk_j+q_j)|y_j|^2
-\sum_{j,k=1}^{\n}\sum\limits_{{\bf
e}=(v_j,\,v_k)\in\cA_F}e^{\,i\lan\t ({\bf e}),\,\vt\ran }\,\bar
y_j\,y_k=\lan H(\vt)y,y\ran.
\]
This, \er{CF3'}, \er{CF3}, \er{CF6} and the minimax principle \er{CF1},
\er{CF2} yield for  $1\leq n\leq\n$:
\[
\lb{CF10'} \l_{n+\n_N-\n}^N\geq\min_{S_{n}\ss \C^{\n}}
\max_{\|y\|=1 \atop {x\in S_{n}}}\lan H(\vt)y,y\ran=\l_n(\vt)\,,
\]
\[
\lb{CF10} \l_n^N\leq\max_{S_{\n-n+1}\ss\C^{\n}} \min_{\|y\|=1 \atop
{x\in S_{\n-n+1}}}\lan H(\vt)y,y\ran=\l_n(\vt)\,.
\]

Second, let $1\leq n\leq\n_D$. Using \er{CF1} and \er{CF2} we write
\[
\lb{CF3.1'} \l_j(\vt)=\min_{S_j\ss \C^{\n}} \max_{\|x\|=1 \atop
{x\in S_j}}\lan H(\vt)x,x\ran\geq\min_{S_j\ss \C^{\n}} \max_{\|x\|=1
\atop {x\in S_j\cap X}}\lan H(\vt)x,x\ran,\qq
j=n+\n-\n_D\,,
\]
\[
\lb{CF3.1} \l_n(\vt)=\max_{S_k\ss \C^{\n}} \min_{\|x\|=1 \atop {x\in
S_k}}\lan H(\vt)x,x\ran\leq \max_{S_k\ss \C^{\n}} \min_{\|x\|=1
\atop {x\in S_k\cap X}}\lan H(\vt)x,x\ran, \qq k=\n-n+1.
\]
For $x\in X$ we have
\[\lb{CF5.1}
\lan H(\vt)x,x\ran=\sum_{j,k=1}^{\n}H_{jk}(\vt)\,\bar
x_j\,x_k=\sum_{j,k=1}^{\n_D}H_{jk}^D\,\bar x_j\,x_k=\lan H_Dx,x\ran,
\]
\[\lb{CF6.1}
\|x\|=\sum_{j=1}^{\n}|x_j|^2=\sum_{j=1}^{\n_D}|x_j|^2.
\]
Then for $1\leq n\leq\n_D$ we may rewrite the inequalities
\er{CF3.1'}, \er{CF3.1} in the form
\[
\lb{CF7.1'} \l_{n+\n-\n_D}(\vt)\geq\min_{S_{n}\ss \C^{\n_D}}
\max_{\|x\|=1\atop {x\in S_{n}}}\lan H_Dx,x\ran=\l_n^D,
\]
\[
\lb{CF7.1} \l_n(\vt)\leq\max_{S_{\n_D-n+1}\ss \C^{\n_D}}
\min_{\|x\|=1 \atop {x\in S_{\n_D-n+1}}}\lan H_Dx,x\ran=\l_n^D.
\]

Combining \er{CF10} and \er{CF7.1} and using \er{snH}, we obtain
\[\lb{FP}
\begin{array}{ll}
\l_n(\vt)\in [\l_n^N,\l_n^D]=J_n,\qqq & n=1,\ldots,\n_D\,, \\[6pt]
\l_n(\vt)\in  [\l_n^N, q_n^\bullet+2\vk_+]=J_n,\qqq & n=\n_D+1,\ldots,\n,
\end{array}
\]
for all $\vt\in\T^d$.

Similarly, from \er{CF10'} and \er{CF7.1'} we obtain
\[\lb{FP1}
\begin{array}{ll}
\l_n(\vt)\in [q_n^\bu, \l_{n+\n_N-\n}^N]=\wt J_n,\qqq & n=1,\ldots,\n-\n_D\,, \\[6pt]
\l_n(\vt)\in [\l_{n+\n_D-\n}^D, \l_{n+\n_N-\n}^N]=\wt J_n,\qqq &
n=\n-\n_D+1,\ldots,\n,
\end{array}
\]
for all $\vt\in\T^d$. The relations \er{FP} and
\er{FP1} prove \er{ell}. \BBox

\

{\bf Proof of Theorem \ref{T2}.}  i) First, we will prove the
estimate \er{est1}. Let $P$ be the projection onto $\ell^2(\partial
V_N)$. Using \er{J1} we have
$$
\begin{aligned}
\sum_{n=1}^{\n}|\s_n(H)|\le
\sum_{n=1}^{\n_D}(\l_n^D-\l_n^N)+\sum_{n=\n_D+1}^{\n}(q_n^\bu+2\vk_+-\l_n^N)\\
=\Tr H_D-\Tr
H_N+\sum_{n=\n+1}^{\n_N}\l_n^N+\sum_{n=\n_D+1}^{\n}(q_n^\bu+2\vk_+)\\
=\sum_{n=\n+1}^{\n_N}\l_n^N+\sum_{n=\n_D+1}^{\n}(q_n^\bu+2\vk_+)-\Tr(P
H_N).
\end{aligned}
$$
Finally, applying \er{H+}, \er{l2.15} to the diagonal entries of $P
H_NP$,   we obtain
\[
\begin{aligned}
\sum_{n=1}^{\n}|\s_n(H)|\le\sum_{n=\n+1}^{\n_N}\l_n^N+\sum_{n=\n_D+1}^{\n}(q_n^\bu+2\vk_+)-
\sum_{n=\n_D+1}^{\n_N}\big(\r_n\,(\vk_n^N-\vk^N_{nn})+q^N_{n}\big)\\
=\sum_{n=\n+1}^{\n_N}\l_n^N+\sum_{n=\n_D+1}^{\n}\big(q_n^\bu+2\vk_+-
\r_n\,(\vk_n-\vk_{nn}+q_{n})\big).
\end{aligned}
\]
Here we have used the following identities
$$
\sum_{n=\n_D+1}^{\n_N}\r_n\,\vk_n^N=\sum_{n=\n_D+1}^{\n}\r_n\sum_{v\in V_n}\,\vk_v^N=\sum_{n=\n_D+1}^{\n}\r_n\vk_n,
$$
$$
\sum_{n=\n_D+1}^{\n_N}\r_n\,\vk_{nn}^N=\sum_{n=\n_D+1}^{\n}\r_n\vk_{nn},\qqq
\sum_{n=\n_D+1}^{\n_N}q_n^N=\sum_{n=\n_D+1}^{\n}\r_nq_n.
$$
Thus, the estimate \er{est1} has proved.

Second, using \er{J1} and \er{wJ1} we have
\begin{multline*}
\sum_{n=1}^{\n}|\s_n(H)|\le
\sum_{n=1}^{\n_D}(\l_n^D-\l_n^N)+\sum_{n=\n_D+1}^{\n-\n_D}(\l_{n+\n_N-\n}^N-
\l_n^N)\\+\sum_{n=\n-\n_D+1}^{\n}(\l_{n+\n_N-\n}^N-
\l_{n-\n+\n_D}^D)=
\sum_{n=\n_D+1}^{\n}\l_{n+\n_N-\n}^N-
\sum_{n=1}^{\n-\n_D}\l_{n}^N\\=\sum_{n=1}^{\n-\n_D}\l_{\n_N-(\n-\n_D)+n}-
\sum_{n=1}^{\n-\n_D}\l_n^N.
\end{multline*}
Thus, the estimate \er{est2} has also proved.

ii) First, we will prove that $0\leq\n_D\leq\n-1$. It is clear that $0\leq\n_D$ and, for example, for the square lattice $\n_D=0$. For the graph shown in Fig.\ref{ff.10} $\n_D=\n-1$. Assume that $\n_D=\n$. Then the graph $\G_N$ contains all bridges of $\G$. This contradicts the construction of $\G_N$.

\setlength{\unitlength}{1.0mm}
\begin{figure}[h]
\centering
\unitlength 1mm 
\linethickness{0.4pt}
\ifx\plotpoint\undefined\newsavebox{\plotpoint}\fi 
\begin{picture}(60,50)(0,0)

\put(10,10){\line(1,0){40.00}}
\put(10,30){\line(1,0){40.00}}
\put(10,50){\line(1,0){40.00}}
\put(10,10){\line(0,1){40.00}}
\put(30,10){\line(0,1){40.00}}
\put(50,10){\line(0,1){40.00}}

\put(10,10){\circle{1}}
\put(30,10){\circle{1}}
\put(50,10){\circle{1}}

\put(10,30){\circle{1}}
\put(30,30){\circle{1}}
\put(50,30){\circle{1}}

\put(10,50){\circle{1}}
\put(30,50){\circle{1}}
\put(50,50){\circle{1}}

\put(10,10){\vector(1,0){20.00}}
\put(10,10){\vector(0,1){20.00}}

\put(7.5,7.5){$\scriptstyle v_\n$}
\put(1,29){$\scriptstyle v_\n+a_2$}
\put(29,7){$\scriptstyle v_\n+a_1$}
\put(20,8){$\scriptstyle a_1$}
\put(6,20){$\scriptstyle a_2$}
\put(14,26.5){$\scriptstyle v_1$}
\put(23.8,17){$\scriptstyle v_{\nu-2}$}

\put(24.0,12){$\scriptstyle v_{\nu-1}$}
\put(20,26.5){$\scriptstyle v_2$}

\put(10,10){\line(2,3){10.00}}
\put(10,10){\line(3,2){15.00}}
\put(10,10){\line(1,3){5.00}}
\put(10,10){\line(3,1){15.00}}
\put(20.5,21){$\ddots$}
\put(15,25){\circle{1}}
\put(20,25){\circle{1}}
\put(25,15){\circle{1}}
\put(25,20){\circle{1}}

\put(30,10){\line(2,3){10.00}}
\put(30,10){\line(3,2){15.00}}
\put(30,10){\line(1,3){5.00}}
\put(30,10){\line(3,1){15.00}}
\put(40.5,21){$\ddots$}
\put(35,25){\circle{1}}
\put(40,25){\circle{1}}
\put(45,15){\circle{1}}
\put(45,20){\circle{1}}

\put(10,30){\line(2,3){10.00}}
\put(10,30){\line(3,2){15.00}}
\put(10,30){\line(1,3){5.00}}
\put(10,30){\line(3,1){15.00}}
\put(20.5,41){$\ddots$}
\put(15,45){\circle{1}}
\put(20,45){\circle{1}}
\put(25,35){\circle{1}}
\put(25,40){\circle{1}}

\put(30,30){\line(2,3){10.00}}
\put(30,30){\line(3,2){15.00}}
\put(30,30){\line(1,3){5.00}}
\put(30,30){\line(3,1){15.00}}
\put(40.5,41){$\ddots$}
\put(35,45){\circle{1}}
\put(40,45){\circle{1}}
\put(45,35){\circle{1}}
\put(45,40){\circle{1}}
\end{picture}
\vspace{-0.5cm} \caption{\footnotesize  $\Z^2$-periodic graph $\G$ with $\n_D=\n-1$.} \label{ff.10}
\end{figure}
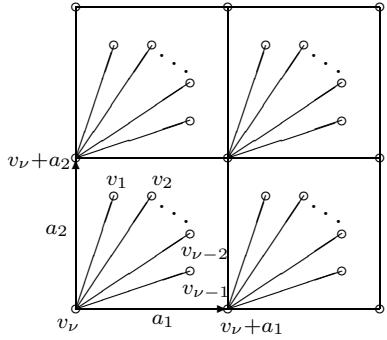

Second, we will show that $\n+d\leq\n_N\leq\n+{\b\/2}$. Due to the identity
$$
V_N=V_F\cup \{u\in V:(u,v)_e\in \cB_N, v\in V_F \},
$$
we have $\n_N=\# V_N\leq\n+{\b\/2}$\,. Since $\Z^d$-periodic graph
$\G$  is connected, the minimal number of vertices from
$V_N\setminus V_F$ is $d$. Thus, $\n+d\leq\n_N$. For example, for
the square lattice $\n_N=\n+{\b\/2}=\n+d$, i.e., the boundaries of
the second inequality in \er{nv} are also achieved. \qq \BBox

\subsection{\lb{S3}Example.}
Consider the Laplacian $H=\D$ on the periodic graph $\G$ shown in
Fig.\ref{ff.0.11}\emph{a}. For each $\vt\in\T^2$ the matrix $\D(\vt)$ defined by (\ref{l2.15'})
has the form
\begin{equation}\label{z2}
\D(\vt)=\left(
\begin{array}{ccc}
   6 & -\D_{12}(\vt)  & -1-e^{-i\vt_2} \\ [2pt]
   -\bar\D_{12}(\vt) & 4 & 0 \\ [2pt]
   -1-e^{i\vt_2} & 0 & 2
\end{array}\right),\qq \D_{12}(\vt)=1+e^{i\vt_1}+e^{-i\vt_2}+e^{i(\vt_1-\vt_2)}.
\end{equation}
The characteristic polynomial of $\D(\vt)$ is given by
$$
\det(\D(\vt)-\l\1_3)=-\l^3+12\l^2+2\,(2c_1c_2+2c_1+3c_2-19)\l
-4(2c_1c_2+2c_1+4c_2-8),
$$
$$
c_1=\cos\vt_1,\qq
c_2=\cos\vt_2.
$$
The spectrum of the Laplacian $\D$ on the periodic graph $\G$
consists of three bands:
\[\lb{spe}
\textstyle\s_1=[0;2],\qqq
\s_2\approx[2{.}5;4],\qqq \s_3\approx[6;9{.}5].
\]

The matrices $H_N$ and $H_D$, defined by \er{H+} -- \er{H-}, in this
case have the form
$$
H_N=\left(
\begin{array}{ccccccc}
6 & -2  & -\sqrt{2} & -2 & -2 & -\sqrt{2} & -2\\
-2 & 4  & 0 & 0 & 0 & 0 & 0 \\
-\sqrt{2} & 0  & 2 & 0 & 0 & 0 & 0\\
-2 & 0  & 0 & 4 & 0 & 0 & 0 \\
-2 & 0 & 0 &  0 & 4 & 0 & 0 \\
-\sqrt{2} & 0 & 0 &  0 & 0 & 2 & 0 \\
-2 & 0 & 0 &  0 & 0 & 0 & 4 \\
\end{array}\right),\qqq H_D=6.
$$
The spectra of the operators $H_N$ and $H_D$ are
$$
\textstyle\s(H_N)\approx \big\{0;
2;2{.}5;4;4;4;9{.}5\big\},\qqq \s(H_D)=
\big\{6\big\}.
$$
Thus, the intervals $J_n$ and $\wt J_n$ defined by \er{J1}, \er{wJ1}
and their intersections $J_n\cap\wt J_n$,  $n\in\N_3$, have the form
$$
\begin{array}{lll}
J_1=[0;6], & \wt J_1=[0,4], \qqq  &
\s_1=[0;2]\ss J_1\cap\wt J_1=
J_1=[0;4],\\[6pt]
J_2=[2;12], \qqq & \wt J_2=[0;4],  &
\s_2\approx[2{.}5;4]\ss J_2\cap\wt J_2=[2;4],\\[6pt]
J_3\approx[2{.}5;12], & \wt J_3\approx[6;9{.}5], &
\s_3\approx[6;9{.}5]=J_3\cap\wt J_3=\wt J_3.
\end{array}
$$

\no \textbf{Remark.} 1) Theorem \ref{T1} determines the existence of the second spectral gap
(see Fig.\ref{ff.0.11}\emph{b}). The intersection of the intervals $J_n$ and $\wt J_n$, $n=1,2,3$, gives more
precise estimates of the spectral band $\s_n(H)$ than one interval
$J_n$. Moreover, for $n=2,3$
the estimate $\s_n(H)\subset J_n$ gives the upper bound
$\l_n(\vt)\leq2\vk_+=12$ that is trivial. But using \er{ell} we obtain more accurate
estimates for the spectral bands. Note that the last spectral band and the last gap are detected precisely, but the first band is estimated too roughly and the first spectral gap is not detected.

2) For the graph shown in Fig.\ref{ff.0.11}\emph{a} the estimates \er{est1}, \er{est2} have the form

\[\lb{est11}
\sum_{n=1}^{3}|\s_n(H)|\le\sum_{n=2}^3(12-h_n)+
\sum_{n=4}^{7}\l_n^N\approx4+12+9{.}5=25{.}5,
\]
\[\lb{est22}
\sum_{n=1}^{3}|\s_n(H)|
\le\sum_{n=1}^{2}\big(\l_{5+n}^N-\l_n^N\big)\approx9{.}5+4-2=11{.}5.
\]
Thus, the estimate \er{est22} is much better than \er{est11}. Finally, we note that  \er{spe} yields
$$
\sum_{n=1}^3|\s_n(H)|\approx(2-0)+(4-2{.}5)+(9{.}5-6)=7{.}0.
$$


\medskip

\footnotesize \no\textbf{Acknowledgments.} \footnotesize Various
parts of this paper were written during Evgeny Korotyaev's stay  in
 Centre for Quantum Geometry of
Moduli spaces (QGM),  Aarhus University, Denmark. He is grateful to
the institute for the hospitality. His study was partly supported by
The Ministry of education and science of Russian Federation, project
07.09.2012 No 8501 and the RFFI grant "Spectral and asymptotic
methods for studying of the differential operators" No 11-01-00458
and the Danish National Research Foundation grant DNRF95 (Centre for
Quantum Geometry of Moduli Spaces - QGM)".

\end{document}